\documentclass[11pt,oneside]{amsart}
 \usepackage{graphicx}
\usepackage{amsfonts, amssymb}
\usepackage{amsmath,amsthm,amscd,xypic}
\usepackage{color}
\usepackage{epic,eepic}
\usepackage{geometry}
\newtheorem{theorem}{Theorem}[section]
\newtheorem{thm}[theorem]{Theorem}
\newtheorem*{thm*}{Theorem}
\newtheorem{lem}[theorem]{Lemma}

\newtheorem{cor}[theorem]{Corollary}
\newtheorem*{cor*}{Corollary}

\newtheorem*{Theorem A}{Theorem A}
\newtheorem*{Theorem B}{Theorem B}
\newtheorem*{Theorem C}{Theorem C}
\newtheorem*{Theorem D}{Theorem D}
\newtheorem*{Theorem E}{Theorem E}

\theoremstyle{definition}

\newtheorem{example}[theorem]{Example}

\newtheorem*{conj*}{Conjecture}

\newtheorem{remark}[theorem]{Remark}

\numberwithin{equation}{section}

\begin{document}

\title{Bounding geometry of loops in Alexandrov spaces}

\author{Nan Li}
\address{Department of Mathematics, University of Notre Dame, Notre Dame, IN 46556}
\email{nli2@nd.edu}
\thanks{Both authors are supported
partially by NSF Grant DMS 0203164 and by a reach found from Capital
normal university.}

\author{Xiaochun Rong}
\address{Department of Mathematics, Capital Normal University, Beijing,
P.R.C.}
\address{Department of Mathematics, Rutgers University, New Brunswick, New Jersey 08854}
\email{rong@math.rutgers.edu}





\begin{abstract}
  For a path in a compact finite dimensional
Alexandrov space $X$ with curv $\ge
\kappa$, the two basic geometric invariants are the length and the
turning angle (which measures the closeness from being a geodesic).
We show that the sum of the two invariants of any loop is bounded
from below in terms of $\kappa$, the dimension, diameter and
Hausdorff measure of $X$. This generalizes a basic estimate of
Cheeger on the length of a closed geodesic in closed Riemannian manifold
([Ch], [GP1,2]). To see that the above result
also generalizes and improves an analogous of the Cheeger type
estimate in Alexandrov geometry in [BGP], we show that for a class of
subsets of $X$, the $n$-dimensional Hausdorff measure and rough
volume are proportional by a constant depending on $n=\dim(X)$.
\end{abstract}

\maketitle

\section*{Introduction}

Let $X$ denote an Alexandrov space with curvature bounded from below,
curv $\ge \kappa$, which is a length metric space such that
each point has a neighborhood in which any geodesic triangle
looks fatter than a comparison triangle in the $2$-dimensional space
form $S^2_\kappa$ of constant curvature $\kappa$.
A motivation for studying
Alexandrov spaces is that the Gromov-Hausdorff limit of a sequence of
Riemannian $n$-manifolds with sectional curvature $\text{sec}\ge \kappa$ is
an Alexandrov space with curv $\ge \kappa$.
A Riemannian manifold with $\text{sec}\ge \kappa$
is an Alexandrov space, but an Alexandrov space in general may
have geometrical or topological singularities.
A basic issue in Alexandrov geometry is to prove results
whose counterparts in Riemannian geometry reply on the
Toponogov triangle comparison theorem ([BGP]).

Let $\gamma: [0,1]\to X$ be a continuous curve. Given a partition,
$P: 0=t_1<\cdots <t_{m+1}=1$ with partition size $|P|=\delta$, let $p_i=\gamma(t_i)$, and let
$\gamma_m$ denote an $m$-broken geodesic i.e., $\gamma_m|_{[t_i,t_{i+1}]}=[p_ip_{i+1}]$ is
a minimal geodesic jointing $p_i$ and $p_{i+1}$. Let
$\theta_i=\pi- \measuredangle p_{i-1}p_ip_{i+1}$. In particular,
$\theta_1=\pi-\measuredangle p_{m+1}p_1p_2$ if $p_{m+1}=p_1$
(the loop case) and $\theta_1=0$ otherwise. Let $\Theta_P(\gamma)=\sum_{i=1}^m\theta_i$.
We define the following number,
$$\Theta(\gamma)=\lim_{\delta\to 0}\sup_{|P|=\delta}\{\Theta_P(\gamma)\},$$
the {\it turning angle} of $\gamma$. For
convenience, we assign $2\pi$ as the turning angle of a trivial
loop. An $m$-broken geodesic $\gamma_m$ has a finite turning angle
$\Theta(\gamma_m)=\sum_{i=1}^m\theta_i$. $\Theta(\gamma)$ measures the closeness of a curve from a geodesic in the following sense. A curve $\gamma$ is a geodesic if and only if $\Theta(\gamma)=0$. If $M$ is a Riemannian
manifold, then for any $C^2$-curve $\gamma\subset M$,
$\Theta(\gamma)=\int^1_0|\nabla _{\gamma'}\gamma'|dt$ is the geodesic curvature (c.f. [AB]). Because a general Alexandrov space
may contain no closed geodesic (nor an $m$-broken
geodesic loop with a small turning angle; e.g., a flat cone), a loop
with the minimal turning angle should be treated as a counterpart
of a closed geodesic on a (closed) Riemannian manifold.

In this paper, $\text{Haus}_n$ will denote the ``normalized''
$n$-dimensional Hausdorff measure such
that $\text{Haus}_n(I^n)=1$, where $I^n$ is the unit $n$-cube in
$\Bbb R^n$. In particular, if $U$ is an open subset of
an $n$-dimensional Riemannian manifold,
$\text{Haus}_n(U)=\text{vol}(U)$. Let $\text{Alex}^n(\kappa)$ be the collection of $n$-dimensional Alexandrov spaces with curvature bounded from below by $\kappa$ and
$$\text{Alex}^n(\kappa,D)=\left\{X\in\text{Alex}^n(\kappa), \;\text{diam}(X)\le D.\right\}$$
The purpose of this paper is to find an explicit upper bound for the volume of $X\in\text{Alex}^n(\kappa, D)$ in terms of $\kappa, D$, $L(\gamma)$ and $\Theta(\gamma)$ for any given loop $\gamma\in X$ (Theorem A or Theorem \ref{Theorem 1.1}).

When $X$
is a closed Riemannian manifold, this generalizes a basic
estimate of Cheeger on the length of a closed geodesic in
[Ch] (see Theorem 0.3), as well as an overlap with a generalization
of Cheeger's basic estimate in [GP1] (1.3 Main Lemma), [GP2] (Lemma 1.5). As an application, we will present
a local injectivity radius estimate (see Theorem B). To see that Theorem A
also generalizes and improves an analogous of the Cheeger type
estimate in Alexandrov geometry ([BGP] Lemma 8.6), we show that for any open subset of $X$, the $n$-dimensional Hausdorff measure and rough
volume are proportional by a constant depending on $n=\dim(X)$.


This implies that
Theorem A generalizes and improves an analogous of
the Cheeger type estimate in [BGP] on the length of an almost
closed geodesic in an Alexandrov space (see Theorem 0.5).

We now begin to state the main results of this paper. A more general form will be proposed in Theorem \ref{Theorem 1.1}.

\begin{Theorem A}
Let $X$ be a complete $n$-dimensional Alexandrov space ($n\ge 2$) with
curv $\ge \kappa$. If $\gamma$ is a loop at $p\in X$ contained in a $r$-ball
$B_r(p)$, then the length and turning angle of $\gamma$ satisfy:
$$L(\gamma)+(n-1)r\cdot \Theta(\gamma)\ge \frac{(n-1)\text{Haus}_n(B_r(p))}
{\text{vol}(S^{n-2}_1)\cdot sn_\kappa^{n-1}(r_0)},$$
where $S^m_1$ denotes an unit $m$-sphere, $r_0=r$
for $\kappa\le 0$ and $r_0=\min\{r,\frac \pi{2\sqrt \kappa}\}$ for
$\kappa>0$, and $sn_\kappa(r)=\frac 1{\sqrt \kappa}\sin \sqrt \kappa
r$, $r$, $\frac 1{\sqrt {-\kappa}}\sinh \sqrt {-\kappa}r$
respectively for $\kappa>0$, $\kappa=0$ and $\kappa<0$.
\end{Theorem A}

The lower bound on the left hand side of the inequality
in Theorem A is optimal in all dimensions; the
inequality becomes an equality when $\gamma$ is a great circle in an
$n$-dimensional spherical $\kappa$-space form and $r=\frac \pi{\sqrt \kappa}$
(note that $\text{vol}(S^n_1)=\frac {2\pi}{n-1}\cdot \text{vol}(S^{n-2}_1)$, $n\ge 2$).
Furthermore, in the case when $X$ contains no closed geodesic, the inequality
is sharp modulo a constant depending only on $n$ (see Example \ref{example1}). Let
$\text{Alex}^n(\kappa,D,v)=\left\{X\in\text{Alex}^n(\kappa,D), \;\text{Haus}_n(X)\ge v\right\}$.

\begin{cor}\label{Corollary 0.1}

Let $X\in \text{Alex}^n(\kappa,D,v)$. For any loop $\gamma$ on $X$,
$$L(\gamma)+\Theta(\gamma)\ge c(n,k,D,v)>0,$$
where $c(n,k,D,v)=\frac {v\cdot
\min\{(n-1),D^{-1}\}}{\text{vol}(S^{n-2}_1)\cdot
sn_\kappa^{n-1}(D_0)}$, $D_0=D$ for $\kappa\le 0$ and
$D_0=\min\{D,\frac \pi{2\sqrt \kappa}\}$ for $\kappa>0$.

\end{cor}

Corollary \ref{Corollary 0.1} reveals a basic geometric property of the loop
space over a compact Alexandrov space $X\in \text{Alex}^n(\kappa,D,v)$: any short loop has
turning angle not small, or equivalently, any loop with small
turning angle is not short.

For $0\le \epsilon<1$, we call a loop, $\gamma$, {\it $\epsilon$-closed geodesic},
if $\Theta(\gamma)\le \epsilon\cdot \frac {v}{D\cdot
\text{vol}(S^{n-2}_1)\cdot sn_\kappa^{n-1}(D_0)}$, where
$D_0=D$ for $\kappa\le 0$ and $D_0=\min\{D,\frac\pi{2\sqrt \kappa}\}$
for $\kappa>0$. A loop $\gamma$
is a closed geodesic if and only if $\gamma$ is $0$-closed geodesic.

For any $\epsilon$-closed geodesic $\gamma$ on $X$, its length can be
bounded from below.

\begin{cor}\label{Corollary 0.2}

Let $X\in\text{Alex}^n(\kappa,D,v)$. If $\gamma$ is a loop $\epsilon$-closed geodesic, then
$$L(\gamma)\ge (1-\epsilon)\cdot \frac {(n-1)v}
{\text{vol}(S^{n-2}_1)\cdot sn_\kappa^{n-1}(D_0)},$$ where
$D_0=D$ for $\kappa\le 0$ and
$D_0=\min\{D,\frac \pi{2\sqrt \kappa}\}$ for
$\kappa>0$.
\end{cor}

We will make a few comments on Theorem A:

(a) In Riemannian geometry, it is often important to bound from
below the length of a closed geodesic. For instance, the following
basic estimate of Cheeger on the length of closed geodesics
plays a crucial role in the classical Cheeger's finiteness
theorem ([Ch]).

\begin{thm}[{Cheeger, [Ch]}]\label{Theorem 0.3}

Let $M$ be a closed $n$-manifold ($n\ge 2$) with sectional curvature
$\text{sec}_M\ge \kappa$ ($\kappa\le 0$) and diameter $D<\infty$. For any closed geodesic $\gamma$,
$$L(\gamma)\ge
\frac{(n-1)\text{vol}(M)}{\text{vol}(S^{n-2}_1)\cdot
sn_\kappa^{n-1}(D)}.$$
\end{thm}

Corollary \ref{Corollary 0.2} reduces to Theorem \ref{Theorem 0.3} when restricting to a closed
geodesic (i.e., $\epsilon=0$) on a Riemannian manifold.

(b) We now state a special case of Theorem A.

\begin{Theorem B}

Let $X\in\text{Alex}^n(\kappa,D,v)$. For any $p,q\in X$ and any minimal geodesics $\gamma_1,\gamma_2$ from $p$ to $q$, the distance between $p$ and $q$ satisfies
$$|pq|\ge \frac{n-1}2\cdot \left[\frac{v}
{\text{vol}(S^{n-2}_1)sn_\kappa^{n-1}
(D_0)}-D\cdot \Theta(\gamma_1*\gamma_2^{-1})\right],$$
where for $\kappa\le 0$ and
$D_0=\min\{D,\frac \pi{2\sqrt \kappa}\}$ for
$\kappa>0$.
\end{Theorem B}

Let $a(n,\kappa,D,v)=\frac{v}
{D\cdot\text{vol}(S^{n-2}_1)\cdot sn_\kappa^{n-1}(D_0)}$.
Observe that if $\Theta(\gamma_1*\gamma_2^{-1}) <a(n,\kappa,D,v)$, then $|pq|\ge c(n,\kappa,D, v)>0$.

With a stronger assumption, Theorem B yields an explicit form comparing 1.3 Main Lemma in [GP1] (c.f. [GP2] Lemma 1.5), which generalizes Theorem \ref{Theorem 0.3}.
Consider a compact Riemannian $n$-manifold $M$. For $p, q\in M$, without loss of generality, let $\gamma_1$ and $\gamma_2$ be two minimal
geodesics from $p$ to $q$ such that $\measuredangle(\overset{\cdot} \gamma_1(1),\overset{\cdot}\gamma_2(1))\ge \measuredangle(\overset{\cdot} \gamma_1(0),\overset{\cdot}\gamma_2(0)) =\pi-2\beta$. Then
$$\Theta(\gamma_1*\gamma_2^{-1})=\measuredangle(\overset{\cdot} \gamma_1(0),-\overset{\cdot}\gamma_2(0))+
\measuredangle(\overset{\cdot} \gamma_1(1),-\overset{\cdot}\gamma_2(1))\le 4\beta.$$
Applying Theorem B, we obtain an explicit lower bound for $|pq|$:

\begin{cor}\label{Corollary 0.4}

Let $M$ be a closed $n$-manifold with $\text{sec}_M\ge \kappa$. Assume
$\max\{\text{diam}(\Gamma_{pq}),
\text{diam}(\Gamma_{qp})\}=\pi-2\beta$, where $0\le\beta<\frac{a(n,\kappa,D,v)}
{4}$. Then
$$|pq|\ge \frac{(n-1)D}{2}\left(\frac{\text{vol}(M)}
{D\cdot\text{vol}(S^{n-2}_1)\cdot sn_\kappa^{n-1}(D_0)}- 4\beta\right)>0,$$
where for $\kappa\le 0$ and
$D_0=\min\{\text{diam}(M),\frac \pi{2\sqrt \kappa}\}$ for
$\kappa>0$.
\end{cor}

Comparing to [GP1] and [GP2], let
$S_p$ be the unit tangent sphere, and let
$\Gamma_{pq}\subseteq S_p$
(resp. $\Gamma_{qp}\subseteq S_q$) denote the subset of vectors tangent to
minimal geodesics from $p$ to $q$ (resp. from $q$ to $p$). For any $\theta>0$,
let $\Gamma_{pq}(\theta)=\{\vec{s}\in S_p, |\vec{s}\,\Gamma_{pq}|_{S_p}<\theta\}$, where $|\vec{s}\,\Gamma_{pq}|_{S_p}$ denotes the distance of $\vec{s}$ to $\Gamma_{pq}$ on ${S_p}$. Then
$\Gamma_{pq}(\frac \pi2+\beta)=S_p$
and $\Gamma_{qp}(\frac\pi2+\beta)=S_q$, by [GP1] 1.3 Main Lemma (c.f. [GP2] Lemma 1.5 for an explicit estimate of $\beta$), $|pq|\ge r(n,\kappa,D)$, where $r(n,\kappa,D)$ is of an implicit form. If there is a closed geodesic through $p$ and $q$, then $\beta=0$, and Corollary \ref{Corollary 0.4} implies
Theorem \ref{Theorem 0.3}.

(c) Theorem A can be useful in analyzing local geometry concerning
the {\it injectivity radius} of a point $p$ ($\text{injrad}_p$) in a complete Riemannian manifold $M$. If $q\in M$ is a nearest cut point to $p$ (consequently, $|pq|=\text{injrad}_p<\infty$), then either $q$ is
a conjugate point to $p$ or there is a geodesic loop $\gamma$
at $p$ passing through $q$. In the later case, $2|pq|=
L(\gamma)$ and $\Theta(\gamma)$ satisfy Theorem A.
In the former case (e.g, no geodesic loop satisfying $L(\gamma)=2|pq|$), a similar estimate can also established
(see Theorem B).

To extend a discussion also including an Alexandrov space $X$, we
introduce the following notions: we call a point $p\in X$ a regular
point, if there is a non-trivial minimal geodesic along any
direction in the space of directions at $p$, $\Sigma_p$. As in the
Riemannian case, we define the {\it cut locus}, $C_p$, at a regular
point as the collection of points $q\in X$ such that $q$ is the furthest
point on a radial curve from $p$ with arc length equal to $|pq|$. Let $q\in C_p$
such that $|pq|=|pC_p|$, which equals to the injectivity radius $\text{injrad}_p$.
Clearly, the gradient-exponential map is a homeomorphism on the ball of
radius $<\text{injrad}_p$. Let $\text{geod}(p,q)=\{[pq]\}$ denote the set
of minimal geodesics, $[pq]$, from $p$ to $q$. We call the following
number in $[0,2\pi]$,
$$\theta_p=\inf_{q\in C_p,\,|pq|=\text{injrad}_p}\{\Theta(\gamma_1*\gamma_2^{-1}),\,\,\gamma_1,\gamma_2\in \text{geod}(p,q)\},$$
the {\it geodesic angle} of $p$. Observe that $\theta_p=0$ if and only
if $2\cdot\text{injrad}_p$ is realized by the length of a closed geodesic
at $p$ and $\theta_p=2\pi$ if and only if there is a unique minimal
geodesic $[pq]$ (When $X$ is a Riemannian manifold,
$\theta_p=2\pi$ implies that $q$
is a conjugate point of $p$.). Hence, $\theta_p$ measures the
existence of such a closed geodesic at $p$.

A consequence of Theorem A is:

\begin{cor}\label{Corollary 0.5}

Let $X$ be a complete $n$-dimensional Alexandrov space ($n\ge 2$) with
curv $\ge \kappa$. If $p\in X$ is a regular point, then for any $r>\text{injrad}_p$,
$$\text{injrad}_p\ge \frac{n-1}2\cdot \left[\frac{\text{Haus}_n(B_r(p))}
{\text{vol}(S^{n-2}_1)sn_\kappa^{n-1}
(r_0)}-r\cdot \theta_p\right],$$
where $r_0=r$ for $\kappa\le 0$ and $r_0=\min \{r,\frac \pi{2\sqrt \kappa}\}$
for $\kappa>0$.
\end{cor}

Corollary \ref{Corollary 0.5} provides a local estimate for $\text{injrad}_p$ in terms
of local geometry when $\theta_p$ is relatively small (e.g.,
$\theta_p <\frac{\text{Haus}_n(B_r(p))}{r\cdot \text{vol}(S^{n-2}_1)
\cdot sn_\kappa^{n-1}(r)}$). On the other hand, $\theta_p$ not
relatively small indicates that geodesics from $p$ to $q$ are
confined in a narrow region.

(d) In [BGP], an analogous of Theorem \ref{Theorem 0.3} in Alexandrov geometry
was obtained, which implies a lower bound on the length of
an {\it almost closed geodesic} i.e., an $m$-broken geodesic loop
$\gamma_m=\{[p_ip_{i+1}]\}_{i=1}^m$
($p_{m+1}=p_1$), with $\Theta(\gamma_m)$ very small while $m$ is
fixed. To state the result, we recall two notions in [BGP]: the
$n$-dimensional rough volume of a subset $K\subseteq X$ is the
limit, $V_{r_n}(K)=\underset{\epsilon\to 0}\lim \epsilon^n\cdot
\beta_X(\epsilon)$, where $\beta_X(\epsilon)=\max\{|\{x_i\}|,\,
\{x_i\}\subseteq K \text{ is an $\epsilon$-discrete net}\}$. Clearly, rough
volume is easier to estimate than the Hausdorff measure and
$\text{Haus}_n(X)\le V_{r_n}(X)$. Consider the following function
in $\kappa$ and $d>0$ defined in [BGP]:
$$\psi(\kappa,d)=\max_{q,p,r\in S^2_\kappa}\left\{\frac{|pr|}{\measuredangle pqr},
\, \, |qp|, |qr|, |pr|\le d,
|pr|\ge 2||qp|-|qr||\right\}.$$

\begin{thm}[{[BGP]}]\label{Theorem 0.6}

Let $X$ be a compact $n$-dimensional Alexandrov space of curv $\ge \kappa$.
If $\gamma_m$ is an $m$-broken geodesic loop, then the $n$-dimensional rough volume,
$$V_{r_n}(X)\le \chi_m(\delta_1,\delta)\cdot d\cdot \psi^{n-1}(\kappa,d),$$
where $d=\text{diam}(X)$,
$$\delta_1=\frac 1{\text{diam}(X)}\max\{|p_ip_{i+1}|,\, 1\le
i\le m\},$$
$\max_i\{\theta_i\}\le
\delta$ and $\chi_m(\delta_1,\delta)$ is a constant depending on $m,
\delta_1$ and $\delta$ such that $\chi_m(\delta_1,\delta)\to 0$ as
$\delta_1, \delta\to 0$ ($m$ fixed).
\end{thm}

Theorem \ref{Theorem 0.6} implies a lower bound on the length of an almost
closed geodesic, implicitly in terms of $n, \kappa, d$ and $V_{r_n}(X)$
(when $m$ fixed and $\delta\to 0$, $\delta_1$ must have a positive
lower bound; see Remark 8.7 in [BGP]). However, because $\chi_m(\delta_1,
\delta)\to \infty$ as $m\to \infty$, Theorem \ref{Theorem 0.6}
fails to imply a lower bound on the length of an $m$-broken geodesic
loop (of length, say one) with $m$ large while $m\delta$ are
very small (so both $\delta_1$ and $\delta$ are small).

In view of the above, it is natural to ask if the sharp estimate
in Theorem A holds in terms of the rough volume. First, the rough
volume is not equivalent to the Hausdorff measure in general. For example, the set of rational numbers in $[0,1]$ has rough volume $1$, while its complement and $[0,1]$ both have rough volume $1$. This also concludes that the rough volume does not have aditivity. However, we
can establish the equivalency for the two measures on the bounded subset which is open or has lower dimensional boundary. Note that this includes the closed set whose Hausdorff measure is zero.
Since we can't find this equivalency in literature, for completeness we give a proof
for the following result.

\begin{Theorem C}
Let $U\subseteq X\in \text{Alex}^n(\kappa)$ be a bounded subset. If $U$ is open or the Hausdorff dimension $\dim_H(\partial U)<n$, then
$$V_{r_n}(U)=c(n)\cdot \text{Haus}_n(U),$$
where $c(n)=\frac{V_{r_n}(I^n)}{\text{Haus}_n(I^n)}=V_{r_n}(I^n)$,
and $I^n$ denotes an Euclidean unit $n$-cube.
\end{Theorem C}

Theorem C can be useful in practice; if one wants to prove a result
involving an estimate for $\text{Haus}_n(X)$, then one reduces to
prove it with $V_{r_n}(X)$, which is much easier to estimate. As
for the value of $c(n)$, except $c(1)=1$ and $c(2)\geq\frac 2{\sqrt
3}$, not much is known.

A consequence of Corollary \ref{Corollary 0.2} and Theorem C is:

\begin{cor}\label{Corollary 0.7}

Let $X$ be a compact $n$-dimensional Alexandrov space ($n\ge 2$) with
curv $\ge \kappa$. If $\gamma$ is an $\epsilon$-closed geodesic, then
$$L(\gamma)\ge (1-\epsilon)\cdot \frac {V_{r_n}(X)}{C(n)\cdot sn_\kappa^{n-1}(D_0)},$$
where $D_0=\text{diam}(X)$ for $\kappa\le 0$ and
$D_0=\min\{\text{diam}(X),\frac \pi{2\sqrt \kappa}\}$ for
$\kappa>0$, and $C(n)=\frac{c(n)\cdot
\text{vol}(S_1^{n-2})}{n-1}$ and $c(n)$ is the constant in
Theorem C.
\end{cor}

Corollary \ref{Corollary 0.7} generalizes and improves Theorem \ref{Theorem 0.6} via providing
an explicit sharp estimate for any $\epsilon$-closed
geodesic (including all $m$-broken geodesic loops with $m\delta$
relatively small).

We conclude the introduction by giving an indication for the
proof of Theorem A. First, it is worth to note that our
arguments also implies a new (metric) proof for Theorem \ref{Theorem 0.3};
which does not require a Riemannian structure. Our approach
is very different from the proof of Theorem \ref{Theorem 0.6} in [BGP]
which follows the lines of the proof of Theorem \ref{Theorem 0.3} in [Ch].
Indeed, we found Theorem A after an unsuccessful attempt to
remove the dependence on $m$ from $\chi_m(\delta_1,\delta)$
in Theorem \ref{Theorem 0.6}.

We take an elementary approach to estimate $\text{Haus}_n(X)$
(in the case that $r=\text{diam}(X)$): expressing $\text{Haus}_n(X)$
as a `Riemann sum', bounding each term and evaluating
the ``Riemann sum'' of the bounds via identifying a proper integrant.
Let $\gamma_m=\{[p_ip_{i+1}]\}_{i=1}^m$ be an $m$-broken geodesic loop
approximating to a loop $c$ in Theorem A, and divide $X=\bigcup_{i=1}^mX_i$
such that $\text{Haus}_n(X)=\sum_{i=1}^m \text{Haus}_n(X_i)$, where
$X_i=\{x\in X\,\,\,
|xp_i|\le |xp_j|, \text{ for all $1\le j\ne i\le m$}\}$. Observe that if
$\gamma_m$ is a closed geodesic and $|p_ip_{i+1}|$ is
sufficiently small, then $X_i$ is like the `union of normal
slices' over $[p_ip_{i+1}]$ (when $X$ is a
Riemannian manifold). So in spirit, we are estimating
$\text{Haus}_n(X)$ via a Riemann sum of a double integral: first
over a normal slice at $\gamma_m(t)$, followed by integral over
$\gamma_m$. To obtain a sharp estimate for $\text{Haus}_n(X_i)$,
we apply a basic Hausdorff measure estimate (see Corollary \ref{Corollary 1.6}),
which bounds the Hausdorff measure of any subset $A\subseteq X$
in terms of the Hausdorff measure of the space of directions
at any point $p\in X$, $|pA|$ and $\text{diam}(A\cup \{p\})$.
The key point in
our proof is an estimate of the upper and lower bound for
$\measuredangle xp_ip_{i+1}-\frac \pi2$, $x\in X_i-\{p_i\}$,
in terms of $|p_ip_{i+1}|, |xp_i|$ and $\theta_i$ (see Lemma \ref{Lemma 1.3}).

The rest of the paper is organized as follows:

In Section 1, we will prove Theorem A.

In Section 2, we will prove Theorem C.

\section{Loops and Hausdorff Measure}

Through out this paper, we will freely use basic notions and properties
(such as the space of directions, rough volume, etc) in Alexandrov
geometry. These can be found in [BGP].

The goal in this section is to prove the following volume estimate which easily implies Theorem A.

\begin{thm}\label{Theorem 1.1}

Let $X\in \text{Alex}^n(\kappa)$ ($n\ge 2$). If $\gamma$ is
a loop at $p$ with $\gamma\subset B_r(p)$, then
$$\text{Haus}_n(B_r(p))\le \text{vol}(S^{n-2}_1)\left [\frac{sn^{n-1}_\kappa
(r_0)}{n-1}L(\gamma)+\Theta(\gamma)\int^r_0 sn_\kappa^{n-1}(t)dt\right],$$
where $r_0=r$ for $\kappa\le 0$ and $r_0=\min\{r,\frac\pi{2\sqrt
\kappa}\}$ for $\kappa>0$.
\end{thm}


%
%

To prove Theorem \ref{Theorem 1.1}, it's sufficient to consider the case that $\gamma$ is a broken geodesic loop. Given an $m$-broken geodesic loop, $p\in \gamma_m=\{[p_ip_{i+1}]\}_{i=1}^m\subset B_r(p)$, let $\theta_i=\pi-\measuredangle p_{i-1}p_ip_{i+1}$, then the turning angle $\Theta(\gamma_m)=\sum_{i=1}^m\theta_i$.
We divide
$B_r(p)$ into $m$ subsets ``centered" at $p_i$,
\begin{align*}
  X_i=\{x\in B_r(p),\,\,\, |xp_i|\le |xp_j|,\text{for all $j\ne i$}\},\qquad 1\le i\le m.
\end{align*}


\begin{center}\begin{picture}(220,150)

  \drawline(0,100)(100,75)
  \drawline(100,75)(200,100)
  \put(100,750){\arc{10}{3.386}{6.04}}
  \qbezier(70,150)(30,87.5)(65,0)
  \qbezier(130,150)(170,87.5)(135,0)
  \qbezier(52.5,40)(100,15)(147.5,40)
  \qbezier(52.5,110)(100,135)(147.5,110)
  \qbezier(57.5,22.5)(100,2.5)(142.5,22.5)
  \qbezier(58,127.5)(100,147.55)(142,127.5)
  \drawline(75,24)(100,75)

  \put(0,90){$p_{i-1}$}
  \put(100,65){$p_{i}$}
  \put(190,90){$p_{i+1}$}
  \put(90,85){$\pi-\theta_i$}
  \put(125,122.5){$A_{j}$}
  \put(121,24){ $A_{j}$}
  \put(77.5,21){ $x$}
  \put(110,46){ $X_i$}
  \put(185,0){ Figure 1}
\end{picture}\end{center}

Clearly, $B_r(p)=\bigcup_iX_i$
and thus $\text{Haus}_n(B_r(p))\le \sum_i \text{Haus}_n(X_i)$. We first introduce a volume estimation formula for certain subsets in an Alexandrov space.

\begin{lem}\label{Lemma 1.2}

Let $B_r(p)\subset X\in \text{Alex}^n(\kappa)$, and
let $[pq]$ denote a geodesic in $X$ from $p$ to $q$. Given
$0\le \alpha\le \pi$, $0\le \theta<\pi$ and $L_1, L_2>0$, for $\eta>0$ arbitrarily small, let
\begin{align*} &A([pq],\alpha,L_1,L_2,\theta)=\{x\in B_r(p)-\{p\},
  \\
  &\qquad\qquad\frac {L_2}{\tan_\kappa|xp|}\le \measuredangle
xpq-\alpha +\frac {36\eta^{\frac 32}}{|\tan_\kappa|xp||^{\frac
32}}\le \frac {L_1}{\tan_\kappa|xp|}+\theta\}.
\end{align*}
Then
\begin{align*}
  &\text{Haus}_n(A)
  \\
  &\quad\le \text{vol}(S_1^{n-2})
  \left[\frac{(L_1+L_2)\text{sn}_\kappa^{n-1}(r_0)}{n-1}
  +\theta\cdot\int_0^r\text{sn}_\kappa^{n-1}(t)dt
  +O(\eta^{\frac 32})\right],
\end{align*}
where $r_0=r$ for $\kappa\le 0$ and $r_0=\min \{r,\frac \pi{2\sqrt \kappa}\}$ for $\kappa>0$.
\end{lem}

In fact, the following lemma shows that each $X_i$ is contained in certain subset shaped as in Lemma \ref{Lemma 1.2}.

\begin{lem}\label{Lemma 1.3}

Let the assumptions be as in Theorem \ref{Theorem 1.1} and $\theta_i$, $X_i$ be defined as in the Figure 1. For $\epsilon>0$,
there is $\eta>0$ such that if $\max_i\{|p_ip_{i+1}|\}<\eta$,
then for any $x\in X_i-\{p_i\}$, the following inequality holds:
\begin{align*}
  -\frac {e^{\epsilon}|p_{i}p_{i+1}|}{2\tan_\kappa|xp_i|}
  -\frac {36\eta^{\frac 32}}{|\tan_\kappa|xp_i||^{\frac 32}}
  &\le \measuredangle xp_ip_{i+1}-\frac \pi2
  \\
  &\le \frac {e^{\epsilon}|p_{i}p_{i-1}|}{2\tan_\kappa|xp_i|}
    +\frac {36\eta^{\frac 32}}{|\tan_\kappa|xp_i||^{\frac 32}}+\theta_i,
\end{align*}
where $\tan_\kappa t=\frac{sn_\kappa t}{sn_\kappa'(t)}$, and when $\kappa>0$ and
$|xp_i|=\frac \pi{2\sqrt \kappa}$, the term $\frac {36\eta^{\frac 32}}
{|\tan_\kappa|xp_i||^{\frac 32}}$ is defined to be zero.
\end{lem}

Assuming Lemma Lemma \ref{Lemma 1.2} and \ref{Lemma 1.3}, we can give a proof for Theorem \ref{Theorem 1.1}.

\begin{proof}[{\bf Proof of Theorem \ref{Theorem 1.1}}]

It's sufficient to prove for an $m$-broken geodesic $\gamma_m$, in which $p$ is one of the vertex. For any $\epsilon>0$, evenly
adding $N(\epsilon)$ `broken' points we may assume that the broken
geodesic $\gamma_m$ satisfies that $|p_ip_{i+1}|<\eta$ for all $i$,
where $\eta$ is given in Lemma \ref{Lemma 1.3}. Put $L^i_1=\frac {e^\epsilon |p_{i-1}p_i|}2$ and
$L^i_2=\frac{e^\epsilon |p_ip_{i+1}|}2$. By Lemma \ref{Lemma 1.3}, we see that $X_i\subseteq A([p_ip_{i+1}],\frac \pi2,
L^i_1, L^i_2,\theta_i)$, and by Lemma \ref{Lemma 1.2},
\begin{align*}
  \text{Haus}_n(X_i)
  &\le \text{Haus}_n(A([p_ip_{i+1}],\frac
  \pi2,L^i_1,L^i_2,\theta_i))
  \\
  &\le \text{vol}(S_1^{n-2})
  \left[\frac{e^\epsilon(|p_{i-1}p_i|+|p_ip_{i+1}|)}2\cdot
  \frac{\text{sn}_\kappa^{n-1}(r_0)}{n-1}\right.
  \\
  &\qquad\qquad\qquad\qquad\qquad
  +\left.\theta_i\cdot\int_0^r
  \text{sn}_\kappa^{n-1}(t)dt
  +O(\eta^{\frac 32})\right].
\end{align*}
Then
\begin{align*} &\text{Haus}_n(B_r(p))
  \le \sum_{i=1}^{m+N(\epsilon)} \text{Haus}_n(X_i)
  \\&\quad\le e^\epsilon\cdot \text{vol}(S_1^{n-2})
  \left[\left(\sum_{i=1}^{m+N(\epsilon)}\frac{|p_ip_{i+1}|+|p_1p_{i-1}|}{2}\right)
  \cdot\frac{\text{sn}_\kappa^{n-1}(r_0)}{n-1}\right.
  \\
  &\quad\quad \left.+\sum_{i=1}^{m+N(\epsilon)}\theta_i\cdot\int_0^r
  \text{sn}_\kappa^{n-1}(t)dt+O(\eta^{\frac 12})\right],\hskip2mm \text{(because
$\eta \approx \frac{L(\gamma_m)}{m+N(\epsilon)}$)}
\end{align*}
and the desired inequality follows when $\epsilon\to 0$, and
thus $N(\epsilon)\to \infty$ and $\eta\to 0$.
\end{proof}

To show Lemma \ref{Lemma 1.2}, we need to divide $A$ (or $X_i$ in our context) into thin annulus $A_j$, and then apply an explicit volume formula for $\kappa$-cones (see Lemma \ref{Lemma 1.4}).

For $\Sigma\in\text{Alex}^{n-1}(1)$, one can construct an $n$-dimensional
Alexandrov space $C_\kappa(\Sigma)$ with curv $\ge \kappa$ (cf. [BGP]):
for $\kappa\le 0$, let $C_\kappa(\Sigma)=(\Sigma\times \Bbb R)/(\Sigma\times
\{0\})$ denote a cone over $\Sigma$, and for $\kappa>0$, let
$C_\kappa(\Sigma)=(\Sigma \times [0,\frac\pi{\sqrt \kappa}])/(\Sigma \times
\{0\}, \Sigma\times \{\frac\pi{\sqrt \kappa}\})$ denote the suspension
over $\Sigma$. We define a metric $d$ on $C_\kappa(\Sigma)$ via the
cosine law in the space form of constant sectional curvature $\kappa$.
For instance, if $\kappa=0$, then for $(x,t), (x',t')\in (\Sigma\times
\Bbb R)/(\Sigma\times \{0\})$,
$$d((x,t),(x',t'))^2=t^2+(t')^2-2tt'\cos |xx'|_{\Sigma}.$$
Note that for any $X\in \text{Alex}^n(\kappa)$ and $p\in X$, the space of
directions $\Sigma_p\in \text{Alex}^{n-1}(1)$, and thus we get
$C_\kappa(\Sigma_p)\in \text{Alex}^n(\kappa)$ for a given $\kappa$. If $k>0$, then
$\text{diam}(C_\kappa(\Sigma))=\frac\pi{\sqrt\kappa}$.

Given $\Sigma\in \text{Alex}^{n-1}(1)$ and $0\le r_1<r_2$, let
$$A_{r_1}^{r_2}(\Gamma)=\{x\in C_\kappa(\Sigma): [px]\in\Gamma \text{ and }
r_1\leq|px|\leq r_2\},$$
where $p$ is the vertex of the $\kappa$-cone $C_\kappa(\Gamma)$ which
is a $\kappa$-suspension for $\kappa>0$ (in particular,
$r_2\le \frac \pi{\sqrt \kappa}$ for $\kappa>0$).

\begin{lem}\label{Lemma 1.4}

Let $A^{r_2}_{r_1}(\Gamma)$ be defined as in the above. Then
$$\text{Haus}_n(A^{r_2}_{r_1}(\Gamma))
  =\text{Haus}_{n-1}(\Gamma)
  \cdot \int^{r_2}_{r_1}sn^{n-1}_\kappa(t)dt.$$
\end{lem}

Lemma \ref{Lemma 1.4} is clear if one assumes the co-area formula for Alexandrov spaces ([BGP],
10.6 in [BBI]). Since we do not find a proof in literature for the co-area formula,
for the completeness we will present an elementary proof using the cosine law
in $\kappa$-space form.

\begin{cor}\label{Corollary 1.5}

$$\text{Haus}_n(B_r(C_\kappa(\Gamma))=
\text{Haus}_{n-1}(\Gamma)\cdot \int^r_0sn^{n-1}_\kappa(t)dt.$$
\end{cor}

\begin{cor}\label{Corollary 1.6}

Let $X\in \text{Alex}^n(\kappa)$. Given any bounded subset $A\subseteq X$,
and $p\in X$, then
\begin{equation}\text{Haus}_n(A)\le \text{Haus}_{n-1}(\Gamma_p(A))\int^{r_2}_{r_1}sn_\kappa^{n-1}
(t)dt,\label{1.2.1}\end{equation}
where $\Gamma_p(A)=\{\uparrow_p^q\in\Sigma_p: q\in A\}$,
$r_1=\min_{x\in A}\{|px|\}$ and $r_2=\max_{x\in A}\{|xp|\}$.
\end{cor}

Corollary \ref{Corollary 1.6} may be viewed as an explicit (Hausdorff measure) version of the
comparison theorem in [BGP] Lemma 8.2. One can also see it from Corollary
10.13 in [BGP] assuming the co-area formula for Alexandrov spaces.

\begin{proof}[{\bf Proof of Lemma \ref{Lemma 1.2}}]

Let $A=A([p,q], \alpha,L_1,L_2,\theta)$. Given a partition for
$[0,1]: 0=a_0<a_1<\cdots <a_N=1$, let $r_j=a_jr$, $A_j=\{x\in
A,\,\,\, r_j\le |xp|\le r_{j+1}\}$, $1\le j\le N$. If $\kappa>0$ and
$d>\frac \pi{2\sqrt \kappa}$, we will chose $\{a_j\}$ such that some
$r_j=\frac \pi{2\sqrt \kappa}$ (note that some $A_j$ may be an empty
set; for instance, if $\theta=0$, then $A_j=\varnothing$ when
$r_j>\frac \pi{2\sqrt \kappa}$ because
$\tan_\kappa|xp_i|<0$).
For $x\in A_j$,
$$-\frac {L_2}{\tan_\kappa|xp|}-\frac {36\eta^{\frac 32}}{|\tan_\kappa|xp|
|^{\frac 32}}\le \measuredangle
xpq-\alpha \le \frac {L_1}{\tan_\kappa|xp|}+\theta+\frac {36\eta^{\frac 32}}
{|\tan_\kappa|xp||^{\frac 32}}$$
implies
\begin{align}
-\frac {L_2}{\tan_\kappa(c_j)}-\frac {36\eta^{\frac 32}}{|\tan_\kappa|c_j|
|^{\frac 32}}
&\le \measuredangle
xpq-\alpha
\notag\\
&\le \frac {L_1}{\tan_\kappa(c_j)}+\theta+\frac
{36\eta^{\frac 32}}{|\tan_\kappa|c_j||^{\frac 32}},
\label{1.4.1}\end{align} where
$c_j=r_{j+1}$ when $\kappa\le 0$ or $\kappa>0$ and $r_{j+1}\le \frac
\pi{2\sqrt \kappa}$, otherwise $c_j=r_j$. Let $\Gamma_j=\{[xp]\in
\Sigma_p(X),\ x\in A_j\}$. Because
$\text{curv}(\Sigma_{[pq]}(\Sigma_p))\ge 1$,
$\text{vol}(\Sigma_{[pq]}(\Gamma_j))\leq \text{vol}(S^{n-2}_1)$,
where $\Sigma_{[pq]}(\Gamma_j)$ denotes the space of directions of
$\Gamma_j$ at $[pq]\in \Gamma_j$. Applying Corollary \ref{Corollary 1.6} to $\Gamma_j$
at $[pq]$, by $\text{curv}(\Sigma_p)\ge 1$ and (\ref{1.4.1}) we have
\begin{align}
\text{Haus}_{n-1}(\Gamma_j)
&\le \text{vol}(\Sigma_{[pq]}(\Gamma_j))\cdot \int^{\alpha_1}_{\alpha_2}\sin^{n-2}(t)dt
\notag \\
&\le
\text{vol}(S_1^{n-2})\cdot\left(\frac{L_1+L_2}{\tan_\kappa(c_j)}
+\theta+\frac {72\eta^{\frac 32}}{|\tan_\kappa(c_j)|^{\frac 32}}\right),
\label{1.4.2}\end{align}
where $\alpha_1=\alpha+\frac
{L_1}{\tan_\kappa(c_j)}+\theta+\frac {36\eta^{\frac
32}}{|\tan_\kappa|c_j||^{\frac 32}}$ and $\alpha_2=\alpha-\frac
{L_2}{\tan_\kappa(c_j)}-\frac {36\eta^{\frac
32}}{|\tan_\kappa|c_j||^{\frac 32}}$.
For $\epsilon>0$, when $\triangle_j=r_{j+1}-r_j$ is sufficiently small, we
may assume that $\frac {sn^{n-1}_\kappa(r_{j+1})}{sn_\kappa(r_j)}\le e^\epsilon
sn^{n-2}_\kappa(r_j)$.

Case 1. Assume $\kappa\le 0$ or $\kappa>0$ and $d\le \frac \pi{2\sqrt \kappa}$.
By applying Corollary \ref{Corollary 1.6} to $A_j$: from (\ref{1.4.2}) we get
\begin{align}
  \text{Haus}_n(A_j)
  &\le \text{Haus}_{n-1}(\Gamma_j)\int_{r_j}^{r_{j+1}}
  \text{sn}_\kappa^{n-1}(t)dt
  \notag\\
  &\le \text{Haus}_{n-1}(\Gamma_j) (r_{j+1}-r_j)
  \text{sn}_\kappa^{n-1}(c_j)
  \notag\\
  &\leq \text{vol}(S_1^{n-2})
  \left(\frac{L_1+L_2}{\tan_\kappa(c_j)}+\theta
  +\frac {72\eta^{\frac 32}}
  {|\tan_\kappa(c_j)|^{\frac 32}}\right)
  \text{sn}_\kappa^{n-1}(c_j)\Delta_j
  \notag\\&
  \leq e^\epsilon\cdot \text{vol}(S_1^{n-2})
  \left[(L_1+L_2)\text{sn}_\kappa^{n-2}(c_j)
  \text{sn}'_\kappa(c_j)+\theta
  \cdot\text{sn}_\kappa^{n-1}(c_j)\right.
  \notag\\
  &\quad \left.+72\eta^{\frac 32}
  \text{sn}_\kappa^{n-\frac 52}(c_j)
  \cdot |\text{sn}'_\kappa(c_j)|^{\frac 32}\right] \Delta_j.
\label{1.4.3}\end{align} Then
\begin{align}
  e^{-\epsilon}\cdot\text{Haus}_n&(A)
  =e^{-\epsilon}\cdot\sum_{j=1}^N\text{Haus}_n(A_j)
  \notag\\ & \le \text{vol}(S^{n-2}_1)(L_1+L_2)
  \sum_{j=0}^N sn_\kappa^{n-2}(c_j)sn_\kappa'(c_j)\Delta_j
  \notag\\
  & + \theta\sum_{j=0}^Nsn_\kappa^{n-1}(c_j)\Delta_j
  +72\eta^{\frac 32}\sum_{j=0}^N\text{sn}_\kappa^{n-\frac 52}(c_j)
  \cdot |\text{sn}'_\kappa(c_j)|^{\frac 32}\Delta_j.
\label{1.4.4}\end{align}
Finally, view (\ref{1.4.4}) as Riemann sum of some integrals and let $N\to
\infty$. Note that for $n=2$, $\int^r_0sn_\kappa^{-\frac 12}(t)\cdot
|sn_\kappa'(t)|^{\frac 32}dt<\infty$ because $sn^{-\frac
12}_\kappa(t)=t^{-\frac 12}+o(t)$, we get
\begin{align*} &\text{Haus}_n(A)
  \le e^\epsilon\cdot \text{vol}(S_1^{n-2})
  \left[(L_1+L_2)
  \int_0^{r_0}\text{sn}_\kappa^{n-2}(t)\text{sn}'_\kappa(t)dt\right.
  \\
  &\qquad\qquad\quad\left.+\theta\cdot\int_0^r\text{sn}_\kappa^{n-1}(t)dt
  +72\eta^{\frac 32}\int_0^r\text{sn}_\kappa^{n-\frac 52}(t)\cdot
  |\text{sn}'_\kappa(t)|^{\frac 32}dt\right]
  \\
  &\quad= \text{vol}(S_1^{n-2})
  \left[e^\epsilon\cdot \frac{(L_1+L_2)\text{sn}_\kappa^{n-1}(r_0)}{n-1}
  +\theta\cdot\int_0^r\text{sn}_\kappa^{n-1}(t)dt
  +O(\eta^{\frac 32})\right]
\end{align*}
Letting $\epsilon \to 0$, we see the desired result.

Case 2. Assume $\kappa>0$ and $d>\frac \pi{2\sqrt \kappa}$. For
$A_j$ with $c_j\le \frac \pi{2\sqrt \kappa}$, the estimate in
(\ref{1.4.3}) is still valid. If $c_{j}>\frac \pi{2\sqrt \kappa}$, then we
modify the estimate (\ref{1.4.2}) by throwing out the negative term with
``$\text{tan}_\kappa(c_{j})\leq 0$'', and obtain
\begin{align}
&\text{Haus}_n(A_j)
\notag\\
&\quad\le e^\epsilon\cdot \text{vol}(S^{n-2}_1)
[\theta \cdot sn^{n-1}_\kappa(c_j)+72\eta^{\frac 32}
\text{sn}_\kappa^{n-\frac 52}(c_j)(\text{sn}'_\kappa(c_j))^2]\triangle_i.
\label{1.4.5}\end{align}
Combining (\ref{1.4.3}) for $c_j\le \frac \pi{2\sqrt \kappa}$ and (\ref{1.4.5}), we derive
\begin{align} \text{Haus}_n(A)&=\sum_{j=1}^NV_{r_n}(A_j)
\notag\\&
\le e^\epsilon\cdot \text{vol}(S^{n-2}_1)(L_1+L_2)
\sum_{j=0}^{r_{j+1}\le \frac \pi{2\sqrt \kappa}} sn_\kappa^{n-2}
(c_j)sn_\kappa'(c_j)\Delta_j
\notag\\
& \hskip4mm +
\theta\sum_{j=0}^Nsn_\kappa^{n-1}(r_j)
\Delta_j+O(\eta^{\frac 32}),\label{1.4.6}\end{align}
In (\ref{1.4.6}), letting $N\to \infty$ and $\epsilon\to 0$, we
get
\begin{align*} &\text{Haus}_n(A)
\\
&\le \text{vol}(S^{n-2}_1)\left [(L_1+L_2)
\int^{r_0}_0 sn_\kappa^{n-2}(t)sn_\kappa'(t)dt+\theta\int^r_0
sn_\kappa^{n-1}(t)dt\right]
\\
&=\text{vol}(S^{n-2}_1)\left
[\frac{(L_1+L_2)sn_\kappa^{n-1}(r_0)}{n-1}+\theta\int^r_0sn_\kappa^{n-1}(t)dt\right].
\end{align*}
\end{proof}

\begin{proof}[{\bf Proof of Lemma \ref{Lemma 1.3}}]

For $\epsilon>0$, we may chose $\eta$ small so that for all $i$,
$\frac{|p_ip_{i+1}|}2<\eta$ implies that
$\tan_\kappa\frac{|p_ip_{i+1}|}2\le e^\epsilon\cdot \frac{|p_ip_{i+1}|}2$.
We first claim that
\begin{align}\cos \tilde \measuredangle xp_ip_{i+1}\le  \frac{e^\epsilon\cdot
|p_ip_{i+1}|}{2\tan_\kappa(|xp_i|)},\label{1.3.1}\end{align}
where $\tilde \measuredangle xp_ip_{i+1}$ denotes the corresponding
angle in the comparison triangle $\tilde \triangle xp_ip_{i+1}
\subset S^2_\kappa$. The proof of the claim relies on the cosine law
in the $\kappa$-space form. We will give a proof for the case $\kappa=0, \kappa=-1$ and $\kappa=1$. The general case follows by an analogue modification.

Case 1. Assume $\kappa=0$. By the cosine law and by
the fact that $|xp_i|\le |xp_{i+1}|$, we derive
\begin{align} \cos \tilde \measuredangle xp_ip_{i+1}
&=\frac{|xp_i|^2+|p_ip_{i+1}|^2- |xp_{i+1}|^2}{2|xp_i|\cdot
|p_ip_{i+1}|}
\notag\\
&\le
\frac{|xp_i|^2+|p_ip_{i+1}|^2-|xp_i|^2}{2|xp_i|\cdot
|p_ip_{i+1}|}=\frac{|p_ip_{i+1}|}{2|xp_i|}=
\frac{|p_ip_{i+1}|}{2\tan_0(|xp_i|)}.
\label{1.3.2}\end{align}

Case 2. Assume $\kappa=-1$. By the cosine law and
$|xp_i|\le |xp_{i+1}|$, we derive
\begin{align}
\cos \tilde \measuredangle xp_i&p_{i+1}
=\frac {\cosh |xp_i|\cosh |p_ip_{i+1}|-\cosh |xp_{i+1}|} {\sinh
|xp_i|\sinh |p_ip_{i+1}|}
\notag\\
&\le \frac {\cosh |xp_i|} {\sinh
|xp_i|}\cdot \frac {\cosh |p_ip_{i+1}|-1}{\sinh
|p_ip_{i+1}|}=\frac{\tanh \frac{|p_ip_{i+1}|}2}{\tanh |xp_i|} \le
\frac{|p_ip_{i+1}|}{2\tanh |xp_i|}.\label{1.3.3}\end{align}

Case 3. Assume $\kappa=1$. Again by the cosine law and
$|xp_i|\le |xp_{i+1}|$, we derive:
\begin{align}
\cos \tilde \measuredangle xp_i&p_{i+1}
=\frac{\cos |xp_{i+1}|-\cos |xp_i|\cos
|p_ip_{i+1}|}{\sin |xp_i|\sin |p_ip_{i+1}|}
\notag\\
&\le \frac{\cos |xp_i|-
\cos |xp_i|\cos |p_ip_{i+1}|}{\sin |xp_i|\sin |p_ip_{i+1}|}
\notag\\
&=\frac{\cos |xp_i|2\sin^2 \frac {|p_ip_{i+1}|}2}
{\sin |xp_i|2\sin \frac{|p_ip_{i+1}|}2\cos \frac{|p_ip_{i+1}|}2}
=\frac{\tan \frac{|p_ip_{i+1}|}2}{\tan |xp_i|}
\le \frac{e^\epsilon\cdot |p_ip_{i+1}|}{2\tan|xp_i|}.
\label{1.3.4}\end{align}

By now, (\ref{1.3.1}) follows from (\ref{1.3.2})--(\ref{1.3.4}). Next, we shall show
that the inequality, $u\geq\cos\alpha$, implies
\begin{align}\alpha\geq \frac\pi2-u-36|u|^{\frac 32}.\label{1.3.5}\end{align}
(this will give the left hand side inequality in Lemma \ref{Lemma 1.3}.) Note
that in our case, we may assume $0\le \alpha\le \pi$. Thus, if
$u\geq 1$ or $u\leq -1$, then (\ref{1.3.5}) holds. On the other hand, for
$u\in (-1,1)$, it's sufficient to show $\cos^{-1}u\geq
\frac\pi2-u-36|u|^{3/2}$, equivalently, the function
$$f(u)=u+36|u|^{3/2}-\frac\pi2+\cos^{-1}u\ge 0.$$
By direct calculation,
$$f'(u)=1+
54\cdot \text{sign}(u)|u|^{1/2}-\frac{1}{\sqrt{1-u^2}},\qquad
f''(u)=\frac{27}{|u|^{1/2}}-\frac{u}{(1-u^2)^{3/2}}.$$
For $-1<u< \frac{5\sqrt {13}-1}{18}$, it's easy to see that $f''(u)>0$ and $u=0$ is the only critical point for $f(u)$. Consequently, $f(0)$ is the global minimum for $0<u<\frac{5\sqrt {13}-1}{18}$. For $\frac{5\sqrt {13}-1}{18}<u<1$, $f''(u)<0$ and thus the minimum of $f(u)$ is achieved at the end points. Note that $f(0)=0$ and $f(1)>0$, we get that $f(u)\geq 0$ for
all $u\in(-1,1)$. Plugging in (\ref{1.3.5}) with $\alpha=\measuredangle
xp_ip_{i+1}$ and $u=\frac{e^\epsilon\cdot
|p_ip_{i+1}|}{2\tan_\kappa|xp_i|}$, we obtain
\begin{align} \measuredangle xp_ip_{i+1}
    &\ge \frac \pi2
-\frac {e^{\epsilon}|p_{i}p_{i+1}|}{2\tan_\kappa|xp_i|}
    -36\left(\frac
    {e^{\epsilon}|p_{i}p_{i+1}|}{2|\tan_\kappa|xp_i||}\right)^{3/2}
    \notag\\
    &\ge \frac \pi2
    -\frac {e^{\epsilon}|p_{i}p_{i+1}|}{2\tan_\kappa|xp_i|}
    -\frac {36\eta^{3/2}}{|\tan_\kappa|xp_i||^{3/2}}.
\label{1.3.6}\end{align}
Similarly applying $|xp_i|\leq|xp_{i-1}|$ to the above 3 cases, we obtain
  \begin{align}\measuredangle xp_ip_{i-1}
    \ge \frac \pi2
    -\frac {e^{\epsilon}|p_{i}p_{i-1}|}{2\tan_\kappa|xp_i|}
    -\frac {36\eta^{3/2}}{|\tan_\kappa|xp_i||^{3/2}}.\label{1.3.7}
  \end{align}
Plugging (\ref{1.3.6}), (\ref{1.3.7}) and $\measuredangle
p_{i-1}p_ip_{i+1}=\pi-\theta_i$ into the condition (B) in [BGP]:
$$\measuredangle p_{i-1}p_ip_{i+1}+\measuredangle xp_ip_{i-1}+\measuredangle
xp_ip_{i+1}\le 2\pi,$$
we get the right hand side of the inequality in Lemma \ref{Lemma 1.3}.
\end{proof}




%
%
%
%

As mentioned in the Introduction (see Theorem \ref{Theorem 0.6} and comments
following it), we did not success in an early attempt to modify the
proof of Theorem \ref{Theorem 0.6} in [BGP] in order to remove the dependence on
$m$ from $\chi_m(\delta_1,\delta)$ and factor out $L(\gamma_m)$
from $\chi_m(\delta_1,\delta)$. We like to conclude this section by
explaining the reason for this failure. The proof in [BGP] is,
following the idea in [Ch], to divide $X$ into two parts and
estimate their rough volumes: one part, $U_{\delta_1}$, is like a
$\delta_1$-tube around $\gamma_m$, and the other part,
$X-U_{\delta_1}$. Since points in $X-U_{\delta_1}$ is a definite
distance away from $\{p_i\}$, this allowed [BGP] to have an
estimate for the diameter of the directions pointing to points
in $X-U_{\delta_1}$, in terms of $\delta_1, \delta$ and $m$.
Unfortunately, the rough volumes of two parts in terms of
$\delta_1$ are in different order, that makes it impossible
to remove the dependence on $m$, nor to factor $L(\gamma_m)$,
from $\chi_m(\delta_1,\delta)$.

\section{Hausdorff Measure and Rough Volume}

Our proof of Theorem C relies on the local structure of an Alexandrov
space, which we briefly recall (see [BGP] for details). The
notion of an $(n,\delta)$-strainer maybe viewed as a counterpart of
a normal coordinate on a Riemannian manifold, defined as follows:
for $p\in X$, $n$-pairs of points $\{(p_i,q_i)\}_{i=1}^n$ is
called an $(n,\delta)$-strainer at $p$, if
$$\measuredangle p_ipp_j-\frac \pi2<\delta,\quad \measuredangle p_ipq_i-
\pi<\delta,\quad \measuredangle q_ipq_j-\frac \pi2<\delta. \quad
(1\le i\ne j\le n)$$
We call the number, $\rho=\min\{|pp_i|,|pq_i|\}$, the radius
of the $(n,\delta)$-strainer. By the continuity, the subset of
points with an $(n,\delta)$-strainer is open in $X$. Let $S_\delta$
denote the set of points admitting no $(n,\delta)$-strainer. Then
$S_\delta$ is a closed subset whose Hausdorff dimension $\dim_H(S_\delta)\le n-1$.

Given a bounded set $U\subseteq X\in \text{Alex}^n(\kappa)$, we divide $U$ into the ``regular'' part $U-S_\delta$ and the ``singular'' part $S_\delta$. On the regular part, we have

\begin{lem}[{[BGP] Theorem 9.4}]\label{Lemma 3.1}

Let $X\in \text{Alex}^n(\kappa)$. If $p\in X$ has an $(n,\delta)$-strainer
with radius $\rho>0$, then there are $\epsilon=\epsilon(n,
\delta,\rho)>0$ and $\eta(n,\delta,\rho)>0$ such that $B_\eta(p)$ is $e^\epsilon$ bi-Lipschitz
to an open subset in $\Bbb R^n$. Moreover, $\epsilon\to 0$ as
$\delta\to 0$.
\end{lem}

For our convenience, we call a subset $U$ a region in a metric space with Hausdorff dimension $n$ if the interior of $U$ is non-empty and $\dim_H(\partial U)<n$. By Lemma \ref{Lemma 3.5}, if $U\subseteq X\in \text{Alex}^n(\kappa)$ is a bounded region, then $V_{r_n}(U)=V_{r_n}(\overset{\circ}U)$. In the following we show that Theorem C is true if $X=\mathbb R^n$, In particular, $U$ has no singular point.

\begin{lem}\label{Lemma 3.6}

Let $U\subset \Bbb R^n$ be a bounded region. Then
$$V_{r_n}(U)=c(n)\cdot \text{Haus}_n(U),$$
where $c(n)=\frac{V_{r_n}(I^n)}{\text{Haus}_n(I^n)}$ and $I^n$ is
a unit $n$-cube in $\Bbb R^n$.
\end{lem}

\begin{proof} By Lemma \ref{Lemma 3.5}, it's sufficient to prove for a bounded open set $U$. Note that $\text{Haus}_n(I^n(r))
=r^n\cdot \text{Haus}_n(I^n)$ and $V_{r_n}(I^n(r))=r^n\cdot
V_{r_n}(I^n)$, and thus for any $r>0$,
\begin{align}
V_{r_n}(I^n(r))=c(n)\cdot\text{Haus}_n(I^n(r))
\label{Lemma 3.6.e1}.
\end{align}
It's clear that
\begin{align}
V_{r_n}(I^n_1(r_1)\cup I^n_2(r_2))
=V_{r_n}(I^n_1(r_1)) +V_{r_n}(I^n_2(r_2))
\label{Lemma 3.6.e2}.
\end{align}

We approximate $U$ by finite union of $n$-cubes, whose interior has no overlap with each other. Let $T_j$ and $T_k$ be such approximation satisfying
\begin{align*}
&T_1\subset T_2\subset\cdots \subset T_j \subset\cdots
U\cdots\subset W_k \subset \cdots\subset W_2\subset W_1
\\
&\text{and}\quad\underset j\cup T_j=U=\overset k\cap W_k.
\end{align*}

By (\ref{Lemma 3.6.e2}),
\begin{align*}
  V_{r_n}(U)
  &\ge V_{r_n}(T_j)=\sum_{\alpha\in T_j} V_{r_n}(I^n_{\alpha})
  \\
  &=\sum_{\alpha\in T_j} c(n)\cdot \text{Haus}_n(I^n_{\alpha})
  =c(n)\cdot \text{Haus}_n(T_j).
\end{align*}
Similarly,
$$V_{r_n}(U)\le c(n)\cdot \text{Haus}_n(W_k).$$
Letting $j,k\to\infty$, we get the desired equality.
%
%
%
%
%
%
\end{proof}

Using Lemma \ref{Lemma 3.1} and \ref{Lemma 3.6}, one can get the equivalence for the regular part in $U$. As mentioned in the introduction, for any set $S$, $\text{Huas}_n(S)=0$ may not imply $V_{r_n}(S)=0$. We shall show that this is true in our context (see Lemma \ref{Lemma 3.5}).

Lemma \ref{Lemma 3.3} will be used to improve the following rough volume estimate and get Corollary \ref{Corollary 3.4}. This corollary will be used to deal with the singular part in $U$ (i.e., show Lemma \ref{Lemma 3.5}). Comparing Corollary \ref{Corollary 3.4} with Corollary 8.4 in [BGP], the latter one has the form $V_{r_n}(B_r(p))\le c(n,\kappa,r)$, which is inadequate in our approach for Lemma \ref{Lemma 3.5}.

\begin{lem}[{[BGP], Lemma 8.2}]\label{Lemma 3.2}

Let $X\in\text{Alex}^n(\kappa)$. Given any subset $A\subseteq X$, and $p\in M$,
$$V_{r_n}(A)\le2d_1\psi^{n-1}(\kappa,d)V_{r_{n-1}}(\Gamma_p),$$
where $d_1=\text{diam}(A\cup \{p\})$, $d=\max_{x\in A}\{|px|\}-
\min_{x\in A} \{|px|\}$ and $\Gamma_p\subseteq \Sigma_p$ consists of
geodesic $[pa]$ for every point $a\in A-\{p\}$.
\end{lem}


\begin{lem}\label{Lemma 3.3}

The function $\psi(\kappa,d)$ satisfies the following inequalities:
$$\frac 23\cdot sn_\kappa(d)\le \psi(\kappa,d)\le 2\cdot sn_\kappa(d),$$ provided
$d<\frac \pi{2\sqrt \kappa}$ when $\kappa>0$, where $sn_\kappa(r)$ is defined in Theorem A.
\end{lem}

We will leave the proof of Lemma \ref{Lemma 3.3} to the end of this section. Combining Lemmas \ref{Lemma 3.2} and \ref{Lemma 3.3}, we get

\begin{cor}\label{Corollary 3.4}

Let $p\in X\in \text{Alex}^n(\kappa)$. Then for any $r>0$,
$V_{r_n}(B_r(p))\le c(n,\kappa)\cdot r^n$, where $c(n,\kappa)>0$ is a
constant depending only on $n$ and $\kappa$.
\end{cor}

\begin{lem}\label{Lemma 3.5}

Let $S\subset X\in \text{Alex}^n(\kappa)$ be a compact subset with $\text{Huas}_n(S)=0$. Then
\begin{enumerate}

  \renewcommand{\labelenumi}{(\ref{Lemma 3.5}.\arabic{enumi})}
  \item $V_{r_n}(S)=0$,
  \item there is a sequence $\mu_i\searrow 0$ such that $V_{r_n}(B_{\mu_i}(S))\to 0$ as $i\to\infty$.
\end{enumerate}
\end{lem}

\begin{proof} We argue by contradiction for (\ref{Lemma 3.5}.1). If not so, then there is a sequence
$\epsilon_i\to 0$, and $\epsilon_i$-net $\{x_i^k\}_{k=1}^{\beta(\epsilon_i)}\subset S$
such that
\begin{align}
\epsilon_i^n\cdot\beta(\epsilon_i)\to V_{r_n}(S)>0.
\label{3.5.e2}
\end{align}

Let $B_j(S)=\{x\in X: \text{there is } h\in S \text{ such that } |xh|<1/j\}$ denote the $j^{-1}$-tubular neighborhood of $S$. Because $S$ is closed, $S\subset \cdots \subset B_2\subset B_1$, and $\bigcap _j B_j=S$. Consequently, \begin{align}\text{Haus}_n(B_j)\to \text{Haus}_n(S)=0.\label{3.5.e1}\end{align}
Given any large $j$,
choose $\epsilon_i\le j^{-1}$, and we have
$$\bigcup_kB_{\frac{\epsilon_i}2}(x_i^k)\subseteq B_j,\qquad
B_{\frac {\epsilon_i}2}(x_i^k)\cap B_{\frac{\epsilon_i}2}(x_i^l)=
\varnothing, \quad k\ne l$$ and thus
\begin{align}
&\beta(\epsilon_i)\cdot \min_k \{\text{Haus}_n(B_{\frac{\epsilon_i}2}
(x_i^k))\}
\notag\\
&\qquad\le \sum_k\text{Haus}_n(B_{\frac{\epsilon_i}2}(x_i^k))
\le \text{Haus}_n(B_j).
\label{3.5.e3}
\end{align}
By Bishop-Gromov relative volume comparison for
Alexandrov spaces ([BGP]), we have that for any $p\in X$ and $r>0$,
$$\text{Haus}_n(B_r(p))\ge \frac {\text{Haus}_n(X)}{\text{vol}
(B^\kappa_{\text{diam}(X)})}\cdot \text{vol}(B^\kappa_r)=c(n,\kappa,X)
\cdot r^n>0.$$
In particular, $\text{Haus}_n(B_{\frac {\epsilon_i}2}(x_i^k))\ge
c(n,\kappa,X)\cdot (\frac {\epsilon_i}2)^n$, and thus (\ref{3.5.e3}) implies
\begin{align}
\text{Huas}_n(B_j)
\ge \beta(\epsilon_i)\cdot c(n,\kappa,X)\cdot (\frac {\epsilon_i}2)^n
=\frac{c(n,\kappa,X)}{2^n}\cdot\epsilon_i^n\beta(\epsilon_i).
\end{align}
Let $\epsilon_i\to 0$, we get a contradiction with (\ref{3.5.e2}) and (\ref{3.5.e1}).

To prove (\ref{Lemma 3.5}.2), by (\ref{Lemma 3.5}.1), we may assume a sequence of $\epsilon_i\to 0$
and a sequence of finite $\epsilon_i$-net $\{x_i^k\}_{i=1}^{\beta(\epsilon_i)}\subset S$ such that
$\epsilon_i^n\cdot \beta(\epsilon_i)\le i^{-1}$. Since $\{B_{\epsilon_i}(x_i^k)\}_{i=1}^{\beta(\epsilon_i)}$
is a finite open cover for $S$, we may assume $0<\mu_i<\epsilon_i$ such that
$$B_{\mu_i}(S)\subseteq \bigcup _kB_{\epsilon_i}(x_i^k),$$
and thus
$$V_{r_n}(B_{\mu_i}(S))\le \sum_k V_{r_n}(B_{\epsilon_i}(x_i^k))\le
\beta(\epsilon_i)\cdot \max_k \{V_{r_n}(B_{\epsilon_i}(x_i^k))\}.$$
By Corollary \ref{Corollary 3.4},
$$V_{r_n}(B_{\epsilon_i}(x_i^k))\le c(n,\kappa)\epsilon_i^n,$$
and thus
$$V_{r_n}(B_{\mu_i}(S))\le c(n,\kappa)\cdot (\epsilon_i^n
\cdot \beta(\epsilon_i))\le i^{-1}\cdot c(n,\kappa).$$
\end{proof}

Since $S_\delta$ is closed and $\dim_H(S_\delta)\le n-1$ for $\delta$ small, by Lemma \ref{Lemma 3.5}, have the following.

\begin{cor} \label{Lemma 3.5.cor}
Let $X\in\text{Alex}^n(\kappa)$. Then for $\delta>0$ small, $V_{r_n}(S_\delta)=0$ and there is a sequence $\mu_i\searrow 0$ such that $V_{r_n}(B_{\mu_i}(S_\delta))\to 0$ as $i\to\infty$.
\end{cor}

Now we are ready to prove Theorem C.

\begin{proof}[{\bf Proof of Theorem C}]

Due to Lemma \ref{Lemma 3.5}, it's sufficient to prove for a bounded open set $U$. Fix small $\delta>0$ and take a sequence $\mu_i\searrow 0$. The idea is to divide $U$ into the disjoint union $B_{\mu_i}(S_\delta)\cup \left(U- B_{\mu_i}(S_\delta)\right)$ and verify that
\begin{align}
  &\lim_{i\to\infty}V_{r_n}(B_{\mu_i}(S_\delta))=0\quad\text{and}
  \label{thmc.e1}\\
  &V_{r_n}(U-B_{\mu_i}(S_\delta))=c(n)\cdot\text{Haus}_n(U-B_{\mu_i}(S_\delta)).
  \label{thmc.e2}
\end{align}
By (\ref{thmc.e1}) and $V_{r_n}(U)\le V_{r_n}(U-B_{\mu_i}(S_\delta))+V_{r_n}(B_{\mu_i}(S_\delta))$, we get
\begin{align*}
V_{r_n}(U)\le \lim_{i\to\infty}V_{r_n}(U-B_{\mu_i}(S_\delta))\le V_{r_n}(U).
\end{align*}
Together with (\ref{thmc.e2}),
\begin{align*}
  V_{r_n}(U)
  &=\lim_{i\to\infty}V_{r_n}(U-B_{\mu_i}(S_\delta))
  \\
  &=\lim_{i\to\infty}c(n)\cdot\text{Haus}_n(U-B_{\mu_i}(S_\delta))
  =c(n)\cdot\text{Haus}_n(U).
\end{align*}

(\ref{thmc.e1}) is satisfied due to Corollary \ref{Lemma 3.5.cor}. It remains to show (\ref{thmc.e2}). For each $\mu_i$, because the closure of $U-B_{\mu_i}(S_\delta)$ is compact, we can conclude that every point in $U-B_{\mu_i}(S_\delta)$ has an
$(n,\delta)$-strainer with radius $\rho=\rho(n,\delta,\mu_i)>0$ (if not,
then there is a sequence $x_j\in U-B_{\mu_i}(S_\delta)$ such that the
$(n,\delta)$-strainer at $x_j$ has radius $\rho_i\to 0$. Passing to
a subsequence, we may assume $x_j\to x\in U-B_{\mu_i}(S_\delta)$.
Because the $(n,\delta)$-strainer at $x$ has radius $\rho>0$,
by definition we see that for large $i$, the $(n,\delta)$-strainer
at $x_j$ has radius at least $\rho/2$, a contradiction). By Lemma \ref{Lemma 3.1}, we may assume that $\eta(\delta,\rho)>0$ and
$\epsilon>0$ such that $B_\eta(p)$ is $e^\epsilon$-bi-Lipschitz embedded to
Euclidean space, and $\epsilon\to 0$ as $\delta\to 0$
and $\eta\to 0$ (equivalently, $\delta\to 0$ and $\mu_i\to 0$).

Now we decompose $U-B_{\mu_i}(S_\delta)$ into countable disjoint small regions:
$U-B_{\mu_i}(S_\delta)=\bigcup_jU_j$, such that each $U_j$ is contained
in an $\frac \eta{10}$-ball.
Let $U_j^e$ be the corresponding subset in $\Bbb R^n$ (or equivalently,
$U^e_i$ denotes an Euclidean metric on $U_j$ which is $e^\epsilon$-bi-Lipschitz
to $U_j$). In particular,
$$e^{-\epsilon}\le \frac{V_{r_n}(U_j)}{V_{r_n}(U_j^e)}\le e^\epsilon,
\qquad e^{-\epsilon}\le \frac{\text{Haus}_n(U_j)}{\text{Haus}_n(U_j^e)}
\le e^\epsilon.$$
Together with Lemma \ref{Lemma 3.6}, we get
$$e^{-2\epsilon}c(n)=e^{-2\epsilon}\cdot \frac{V_{r_n}(U_j^e)}{\text{Haus}_n
(U_j^e)}\le \frac{V_{r_n}(U_j)}{\text{Haus}_n(U_j)}\le e^{2\epsilon}
\frac{V_{r_n}(U_j^e)}{\text{Haus}_n(U_j^e)}=e^{2\epsilon}c(n).$$
Because $V_{r_n}$ is finitely additive, we obtain
$$e^{-2\epsilon}c(n)\sum_j\text{Haus}_n(U_j)\le \sum_jV_{r_n}(U_j)\le e^{2
\epsilon}c(n)\sum_j \text{Haus}_n(U_j),$$
and thus
\begin{align}
e^{-2\epsilon}c(n)\cdot \text{Haus}_n(B_{\mu_i}(S_\delta))
&\le V_{r_n}
(U-B_{\mu_i}(S_\delta))
\notag\\
&\le e^{2\epsilon}c(n)\cdot \text{Haus}_n (U-B_{\mu_i}(S_\delta)).\label{3.7.2}
\end{align}
In (\ref{3.7.2}), letting $\delta\to 0$ and ${\mu_i}\to 0$ (thus $\epsilon\to 0$), we get (\ref{thmc.e2}).
\end{proof}

\begin{remark} \label{remark.sec2.1}
We see that both (\ref{thmc.e1}) and (\ref{thmc.e2}) are verified relying on the Alexandrov structure.
\end{remark}

\begin{proof}[{\bf Proof of Lemma \ref{Lemma 3.3}}]

We will first reduce the proof to the case
when $|qp|=|qr|$ (see
(\ref{3.3.1}) below). We may assume
that $|qp|\geq |qr|$, and let $s$ be a point on the geodesic
from $q$ to $p$ such that $|qs|=|qr|=x$. From the condition that
$2(|qp|-|qr|)\le |pr|$, we derive
$$|pr|-|rs|\leq |ps|=|qp|-|qr|\leq\frac 12|pr|,$$
and thus $|pr|\leq 2|rs|$. From
$$|rs|\leq|pr|+|ps|=|pr|+|qp|-|qr|
  \leq |pr|+\frac 12|pr|,$$
we get that $|pr|\geq\frac 23|rs|$, and therefore
$$\frac{2}{3}\frac{|rs|}\theta\leq \frac{|pr|}\theta \leq
  2\frac{|rs|}\theta, $$
where $\theta=\measuredangle pqr$. In the above inequality,
taking maximum over $p, q, r\in S^2_\kappa$ under
the conditions for $\psi(\kappa,d)$, we get
\begin{align}
\frac 23\max_{q,r,s\in S^2_\kappa}\left\{\frac {|rs|}\theta,\,\, |qs|=|qr|\le
d\right\}&\le \psi(\kappa,d)
\notag\\
&\le 2\max_{q,r,s\in S^2_\kappa}\left\{\frac{|rs|}\theta,\,\,
|qr|=|qs|\le d\right\}.\label{3.3.1}
\end{align}
We claim that for each fixed $x$,
\begin{align}
\max\limits_{|rs|}\left\{\frac{|rs|}
\theta,\;|qr|=|qs|=x\right\}=sn_\kappa x.\label{3.3.2}
\end{align}
Clearly, Lemma \ref{Lemma 3.3} follows from (\ref{3.3.1}) and (\ref{3.3.2}).
In the rest of the proof, we will verify (\ref{3.3.2}).

Case 1. For $k<0$, applying the cosine law to the triangle $\triangle qrs$ we derive
\begin{align*}
  \cosh(\sqrt{-\kappa}|rs|)
  &=\cosh^2(\sqrt{-\kappa}x)-\sinh^2(\sqrt{-\kappa}x)\cos \theta
  \\&=1+\sinh^2(\sqrt{-\kappa}x)(1-\cos \theta)
  \\&=1+2\sinh^2(\sqrt{-\kappa}x)\sin^2{\frac \theta 2},
\end{align*}
and thus
\begin{equation}
  \sinh\frac{\sqrt{-\kappa}|rs|}{2}=
  \sin{\frac \theta 2}\sinh(\sqrt{-\kappa}x).
\label{3.3.3}\end{equation}
Since $\sin z\le z$ and $z\le \sinh z$ for $z>0$, from (\ref{3.3.3}) we
get
$$\frac{\sqrt{-\kappa}|rs|}{2}\leq \sinh\frac{\sqrt{-\kappa}|rs|}{2}=
  \sin{\frac\theta 2}\sinh(\sqrt{-\kappa}x) \leq
  \frac \theta 2\sinh(\sqrt{-\kappa}x),
$$ and thus
$$\frac{|rs|}\theta\leq
  \frac{\sinh(\sqrt{-\kappa}x)}{\sqrt{-\kappa}}.
$$
On the other hand, $|rs|\to 0 \Leftrightarrow \theta\rightarrow 0$. Using (\ref{3.3.3}), we derive
$$\lim_{\theta \rightarrow 0}\frac{|rs|}\theta
  =\lim_{\theta\rightarrow 0}\frac{|rs|}{\sinh\frac{\sqrt{-\kappa}|rs|}{2}}\cdot
  \frac{\sin{\frac \theta 2}\sinh(\sqrt{-\kappa}x)}\theta
  =\frac{\sinh(\sqrt{-\kappa}x)}{\sqrt{-\kappa}}.$$
By now, we can conclude (\ref{3.3.2}) for $k<0$.

Case 2. For $k=0$, applying the cosine law to $\triangle qrs$, we
get that $|rs|=2x\sin \frac \theta 2 \le \theta x$ and thus $\frac
{|rs|}\theta\le x$. On the other hand,
$$\lim_{\theta\to 0}\frac{|rs|}\theta=\lim_{\theta\to 0}\frac {2x\sin
\frac \theta 2}\theta=x.$$
Similarly, we can conclude (\ref{3.3.2}) for $k=0$.

Case 3. For $\kappa>0$, applying the cosine law to $\triangle qrs$, we
get
\begin{align}
\sin\frac{\sqrt{k}|rs|}{2}=
  \sin{\frac\theta 2}\sin(\sqrt{k}x).
\label{3.3.4}\end{align}
By (\ref{3.3.4}), we get
\begin{align} \frac{|rs|}\theta &=\frac {\frac{\sqrt \kappa|rs|}2}{\sin
\frac{\sqrt \kappa|rs|}2}\cdot \frac{\sin \frac {\sqrt \kappa |rs|}2}
{\sqrt \kappa\frac \theta 2}=
\frac {\frac{\sqrt \kappa|rs|}2}{\sin \frac{\sqrt \kappa|rs|}2}\cdot \frac{\sin
\frac \theta 2}{\frac \theta 2}\cdot \frac{\sin (\sqrt \kappa x)}{\sqrt \kappa}.
\label{3.3.5}\end{align}
We claim that
$$\frac {\frac{\sqrt \kappa|rs|}2}{\sin \frac{\sqrt \kappa|rs|}2}\cdot
\frac{\sin \frac\theta 2}{\frac \theta 2}\le 1. $$
Because $\theta\to 0$ if and only if $|rs|\to 0$,
$$\lim_{\theta\to 0}\frac {\frac{\sqrt \kappa|rs|}2}{\sin \frac{\sqrt
\kappa|rs|}2}\cdot \frac{\sin \frac
\theta 2}{\frac \theta 2}=1,$$
and consequently we conclude from (\ref{3.3.5}) that (\ref{3.3.2}) holds for $\kappa>0$.

To see the claim, let $\lambda=\sin (\sqrt \kappa x)$, and rewrite (\ref{3.3.4}) as
$$\sin \frac {\sqrt \kappa|rs|}2=\lambda\sin \frac \theta2,\qquad \qquad
\frac {\sqrt \kappa|rs|}2=\sin^{-1}(\lambda\sin \frac \theta2).$$
Then
$$\frac {\frac{\sqrt \kappa|rs|}2}{\sin \frac{\sqrt \kappa|rs|}2}\cdot
\frac{\sin \frac \theta 2}{\frac \theta 2}=\frac {\sin^{-1}(\lambda\sin
\frac \theta2)}{\lambda\sin \frac \theta2}\cdot\frac {\sin \frac \theta 2}
{\frac \theta2}=\frac{\sin^{-1} (\lambda\sin \frac \theta2)}{\lambda\frac
\theta2}\le 1,$$
because for all $0<\lambda\le 1$ and $0\le \frac \theta2\le \frac \pi2$,
$\lambda\sin \frac \theta2\le \sin (\lambda\frac \theta2)$.
\end{proof}

\begin{example}\label{example1}

We will calculate an example showing that when $X$ contains neither a closed geodesic nor
an almost closed geodesic, the inequality in Theorem A is sharp up to a constant
depending only on $n$.

Consider a sector of angle $\theta$ ($0<\theta<\pi$) in a flat
$2$-disk of radius $d$. We obtain a flat cone, $X^2$, by identifying
the two sides of the sector. Then $\text{vol}(X^2)=\frac 12\theta d^2$.
Let $c$ denote a geodesic loop at a point near the vertex. Then $L(c)<<1$
and $\Theta(c)=\theta$. In this case, the inequality in
Theorem A reads:
$$L(c)+\Theta(c)\cdot d\ge \frac{(2-1)\cdot \text{vol}(X^2)}{\text{vol}(S^0_1)
\cdot d}=\frac \theta 2\cdot d.$$

Let $B^m_d$ denote a closed ball of radius $d$ in $\Bbb R^m$, and let
$X^{m+1}=X^2\times B_d^m$ be the metric product. Then $X^{m+2}$ is
compact Alexandrov space of cur $\ge 0$, and
\begin{align*}
&\text{diam}(X^{m+2})=\sqrt 2d,
\\
&\text{vol}(X^{m+2})
=\text{vol}(X^2)\cdot \text{vol}(B^m_d)=\frac {\text{vol}(S^{m-1}_1)}{2(m+1)}
\cdot \theta\cdot d^{m+2}.
\end{align*}
Let $(p_i,x)\in X^{m+2}=X^2\times B^m_d$ such that $p_i$

converges to the vertex of $X^2$, and let $\gamma_i\subset X^2$
be a sequence of geodesic loops at $p_i$. Then $(\gamma_i,x)\subset X^{m+2}$
is a sequence of geodesic loops such that $L(\gamma_i,x)
=L(\gamma_i)\to 0$ and $\Theta((\gamma_i,0))\equiv \theta$.
Applying Theorem A to $(\gamma_i,0)$ and taking limit as
$i\to \infty$, one gets (we also assume $m=2s$ is even)
\begin{align*} \theta\cdot d&\ge \frac{(m+1)\cdot \text{vol}(X^{m+2})}{(m-1)\cdot
\text{vol}(S^m_1)\cdot d^{m+1}}
\\
&=\frac{\text{vol}(S^{m-1}_1)}
{2(m-1)\cdot \text{vol}(S^m_1)}\cdot \theta\cdot d
\\
&=\frac{\frac{2^{\frac m2}
\pi^{\frac{m-2}2}}{(m-1)!!}}{(m-1)\cdot \frac{\pi^{\frac m2}}{(\frac m2)!}}
\cdot \theta\cdot d
\\
&=\frac 1\pi \cdot \frac 1{2s-1}\cdot \left[\frac{(2s)
\cdot (2s-2)\cdots 4\cdot 2}{(2s-1)\cdot (2s-3)\cdots 3\cdot 1}\right]\cdot
\theta\cdot d
\\
&\ge \frac 1{\pi(2s-1)}\cdot \theta \cdot d.
\end{align*}
\end{example}

\section{Appendix}

In this section, we will give proofs for Lemmas \ref{Lemma 1.4}. The main
ingredient in the proof is the cosine law in the $\kappa$-space form.

\begin{proof}[{\bf Proof of Lemma \ref{Lemma 1.4}}]

Note that for $\kappa>0$, $C_\kappa(\Gamma)$ is a $\kappa$-suspension
over $\Gamma$. If $r_1\ge \frac \pi{2\sqrt \kappa}$, by the symmetry
we see that $\text{Haus}_n(A^{r_2}_{r_1}(\Gamma))=\text{Haus}_n(A^{\frac
\pi{\sqrt \kappa}-r_1}_{\frac \pi{\sqrt \kappa}-r_2}(\Gamma))$. If $r_1
<\frac \pi{2\sqrt \kappa}<r_2$, then similarly we may identify
$$\text{Haus}_n(A^{r_2}_{r_1}(\Gamma))=\text{Haus}_n(A^{\frac\pi{2\sqrt
\kappa}}_{r_1}(\Gamma))+\text{Haus}_n(A^{\frac \pi{\sqrt \kappa}}_{\frac
\pi{\sqrt \kappa}-r_2}(\Gamma)).$$
Hence, without loss of generality we may assume that $r_2\le \frac\pi{2\sqrt \kappa}$.

We will divide $A_{r_1}^{r_2}(\Gamma)$ into small annulus and
express $\text{Haus}_n(A_{r_1}^{r_2}(\Gamma))$ as a Riemannian sum
of the Hausdorff measure of these small annulus. The key in the proof
is an estimate the Hausdorff measure of a small annulus in terms
of the Hausdorff measure of a cross section and the width of
the small annulus (one may view this as a local co-area formula
estimate).

Let $\{t_i\}$ be an $N$-partition of $[r_1,r_2]$ and $\Delta
t=\frac{r_2-r_1}{N}$ be sufficiently small. By the above assumption,
$sn_\kappa(t)$ is increasing in each $[t_1, t_{i+1}]$. Let
$S_t=\{x\in A: |px|=t\}$ and $A_{t_i}^{t_{i+1}}=\{x\in A:
t_i\leq|px|\leq t_{i+1}\}$. Define the product metric
$|(a,u),(b,v)|=\sqrt{|a,b|^2+|u,v|^2}$ over $S_{t_i}\times
[t_i,t_{i+1}]$. Because $S_{t_i}$ is an Alexandrov space and the
normalized $\text{Haus}_n$ has countable additivity, we have
\begin{align}\frac{\text{Haus}_n(S_{t_i}\times [t_i,t_{i+1}])}
  {\text{Haus}_{n-1}(S_{t_i})\cdot(t_{i+1}-t_i)}
  =\frac{\text{Haus}_n(I^n)}{\text{Haus}_{n-1}(I^{n-1})
  \cdot \text{Haus}_{1}(I^{1})}
  =1.\label{2.1.1}
\end{align}
Consider the map $f: A_{t_i}^{t_{i+1}}\rightarrow S_{t_i}\times
[r_1,r_2]$ defined as the following: for $x\in A_{t_i}^{t_{i+1}}$,
let $x'\in S_{t_i}$ be the point on geodesic $[px]$ such that
$|px'|=t_i$, then $f(x)=(x',|px|)$ and
$|f(x_1)f(x_2)|^2=|x_1'x_2'|^2+(|px_1|-|px_2|)^2$.

For any $x_1, x_2\in A_{t_i}^{t_{i+1}}$ Assume $|px_2|\geq|px_1|$.
We will show that
\begin{equation}\frac{|x_1x_2|}{|f(x_1)f(x_2)|}= 1+O(\Delta t)\label{2.1.2}\end{equation}

Applying the following version of cosine law (which can
be easily derived) to the triangle $\triangle px_1x_2$ and $\triangle px'_1x'_2$,
we get that
\begin{align*}
  \text{sn}_\kappa^2\frac{|x_1x_2|}{2}
  &=\text{sn}_\kappa^2\frac{|px_1|-|px_2|}{2}
  +\sin^2\frac{\measuredangle x_1px_2}{2}
  \cdot\text{sn}_\kappa|px_1|\text{sn}_\kappa|px_2|
  \\
  \text{sn}_\kappa^2\frac{|x'_1x'_2|}{2}
  &=\sin^2\frac{\measuredangle x'_1px'_2}{2}
  \cdot\text{sn}^2_\kappa (t_i)
\end{align*}
Since $\measuredangle x_1px_2=\measuredangle x'_1px'_2$,
\begin{align*}\text{sn}_\kappa^2\frac{|x_1x_2|}{2}
  &=\text{sn}_\kappa^2\frac{|px_1|-|px_2|}{2}
  +\frac{\text{sn}_\kappa|px_1|\text{sn}_\kappa|px_2|}{\text{sn}^2_\kappa(t_i)}
  \text{sn}_\kappa^2\frac{|x_1'x_2'|}{2}
  \\
  &=\text{sn}_\kappa^2\frac{|px_1|-|px_2|}{2}
  +(1+O(\Delta t))
  \text{sn}_\kappa^2\frac{|x_1'x_2'|}{2}.
\end{align*}
By the Taylor expansion of
$(\text{sn}^{-1}_\kappa(\sqrt{\text{sn}^2_\kappa(x)+(1+O(\Delta
t))\text{sn}^2_\kappa(y)}))^2$, we get that
\begin{align*}
  |x_1x_2|^2
  &= (|px_1|-|px_2|)^2+|x_1'x_2'|^2+O(\Delta t)|x_1'x_2'|^2
  \\
  &=|f(x_1)f(x_2)|^2+O(\Delta t)|x'_1x'_2|^2.
\end{align*}
which leads to (\ref{2.1.2}). By the cosine law, it's easy to see that
$$\text{Haus}_{n-1}(S_{t_{i}})=sn_\kappa^{n-1}(t_i)\text{Haus}_{n-1}(\Gamma_p).$$
Together with (\ref{2.1.1}) and (\ref{2.1.2}),
\begin{align*} \text{Haus}_n(A_{t_i}^{t_{i+1}})
  &=(1+O(\Delta t))^n \text{Haus}_n(S_{t_i}\times
  [r_1,r_2])
  \\
  &=(1+O(\Delta t))^n \text{Haus}_{n-1}(S_{t_i})\Delta t
  \\
  &=(1+O(\Delta t))^n\text{Haus}_{n-1}(\Gamma_p) sn_\kappa^{n-1}(t_i)\Delta t.
\end{align*}

Summing up the above for $i=0,1,\cdots,N-1$ and let $\max\{\Delta
t\}\rightarrow 0$ we get Lemma \ref{Lemma 1.4}.
\end{proof}

%

%

\vskip 30mm

\bibliographystyle{amsalpha}


\end{document}